\documentclass[12pt,a4paper]{amsart}
\usepackage{amsmath,amssymb,amscd}
\newtheorem{theo}{Theorem}[section]

\newtheorem{defi}[theo]{Definition}
\newtheorem{lem}[theo]{Lemma}
\newtheorem{prop}[theo]{Proposition}

\newtheorem{coro}[theo]{Corollary}
\newtheorem{sublem}[theo]{Sublemma}

\begin{document} 

\begin{abstract}
The aim of this note is to give a simplified proof of the surjectivity of the natural Milnor--Chow homomorphism
$\rho: K^M_n(A) \to CH^n(A,n)$ between Milnor $K$--theory and higher Chow groups for essentially smooth (semi--)local $k$--algebras 
$A$ with $k$ infinite. It implies the exactness of the Gersten resolution for Milnor $K$--theory at the generic point. 
Our method uses the Bloch-Levine moving technique and some properties of the Milnor $K$--theory norm for fields.
\end{abstract}

\title{The Milnor--Chow homomorphism revisited}

\author{Moritz Kerz, Stefan M\"uller-Stach}
\email{m.kerz@dpmms.cam.ac.uk, mueller-stach@uni-mainz.de}
\address{University of Cambridge,
Department of Pure Mathematics and Mathematical Statistics, 
Wilberforce Road, 
Cambridge CB3 0WB,
United Kingdom
}
\address{Johannes Gutenberg--Universit\"at Mainz, Institut f\"ur Mathematik, Staudingerweg 9, 55099 Mainz, Germany} 

\thanks{Supported by Studienstiftung des deutschen Volkes and Deutsche Forschungsgemeinschaft} 
\date{9.1.06}

\maketitle

\section*{Introduction} 

In \cite{G,PEVSMS} the surjectivity of the natural Milnor--Chow homomorphism
$$
\rho: K^M_n(A) \to CH^n(A,n)
$$
between Milnor $K$--theory and higher Chow groups for any essentially smooth (semi--)local $k$--algebra $A$ 
with $k$ infinite was shown. This morphism associates to a symbol
$\{f_1,\ldots,f_n\}$ the graph cycle of the map $f=(f_1,\ldots,f_n)$. \\
In this note we want to give a very simple argument which uses two basic ingredients. 
The first is a new argument derived from fairly elementary properties 
of the norm--map for the Milnor $K$--theory of rings which were sketched in \cite{Ke} 
and build up on the theory of Bass and Tate \cite{BT} (see section 2). The idea in \cite{Ke} was to use 
a Milnor $K$--group which is not induced directly from a ring (or algebra) but only from certain generic elements
of a ring. The same technique can also be used to show the Gersten conjecture for Milnor $K$--theory
of regular (semi--)local rings \cite{Ke}. 
The second input is a standard application of the easy moving lemma of Bloch--Levine \cite{B1, L} which 
implies that we can restrict to the case of cycles with smooth components. This was also used in the proof in \cite{PEVSMS}.

Our main theorem is: 

\begin{theo} \label{main-theorem} Let $A$ be an essentially smooth (semi--)local $k$--algebra with infinite residue fields. 
Then the homomorphism $\rho : K^M_n(A)\to CH^n(A,n)$ is surjective for $n \ge 1$.
\end{theo}

Here for a field $k$ we say that a $k$--algebra $A$ is essentially smooth
if $A$ is the localization of a smooth affine $k$--algebra. 
In fact under the conditions of the theorem one can show $\rho$ is bijective \cite{Ke}.
This theorem has a few beautiful applications: 

\begin{coro} Let $A$ be as above and $X=Spec(A)$ integral (i.e., $A$ a domain with quotient field $F$). 
Then the Gersten resolution for Milnor $K$--theory is exact at the generic point: 
$$
K_n^M(A) {\buildrel i_* \over \to} K_n^M(F) {\buildrel T \over \to} \bigoplus_{x \in X^{(1)}} K_{n-1}^M(k(x)) \to \ldots, 
$$
i.e., $\mathrm{ker}(T)=\mathrm{im}(i_*)$, where $T$ is the tame symbol. 
\end{coro}

The exactness of this complex is well known in codimensions $p \ge 1$ by the work of Gabber and Rost \cite{R} 
and follows in degree zero with the same proof as in \cite{PEVSMS} by comparing with the corresponding sequence for higher Chow groups \cite{B1}. 
Note that the work of Kerz \cite{Ke} implies also the Gersten conjecture, i.e.~the injectivity of $i_*$ for such $k$--algebras $A$. 
There is another nice application to \'etale cohomology: 

\begin{coro}
Assume the Bloch--Kato conjecture \cite{Vn}. Let $A$ be a (semi--)local ring containing an 
infinite field and $l>0$ prime to $char(A)$. Then the graded ring 
$H^*_{\rm et}(A,\mu_l^{\otimes *})$ is generated by elements of degree one.
\end{coro}

{\em Proof.}  First assume that $A$ is essentially smooth over an infinite field. 
The Bloch--Kato conjecture implies that we have an isomorphism 
$$CH^n(A,n)/l {\buildrel \cong \over \to} H^n_{\rm et}(A,\mu_l^{\otimes n})$$ for any $l$ prime to $char(A)$.
Composing with $\rho$ we get a surjective ring homomorphism $K_n^N(A)/l \to H^n_{\rm et}(A,\mu_l^{\otimes n})$ which
shows the corollary in this case, because Milnor $K$--theory is generated in degree one. 

Let $A$ be arbitrary. By a direct limit argument we can assume $A$ to be a localization of an affine algebra. 
Now Hoobler's trick \cite{H} can be applied: There is a Henselian pair $(A',I)$ with $A'$ essentially smooth 
and $A=A'/I$. In this situation $H^n_{\rm et}(A,\mu_l^{\otimes n})$ and $H^n_{\rm et}(A',\mu_l^{\otimes n})$
are isomorphic \cite{H}. The commutative diagram 
$$
\begin{CD}
K^M_n(A')/l @>{nat}>> K^M_n(A)/l\\
\rho @VVV \rho @VVV \\
H^n_{\rm et}(A',\mu_l^{\otimes n}) @>>{nat}> H^n_{\rm et}(A,\mu_l^{\otimes n})
\end{CD}
$$
implies immediately that $\rho:  K_n^N(A)/l \to H^n_{\rm et}(A,\mu_l^{\otimes n})$ is surjective.
\hfill $\Box$

\begin{coro}[Bloch] Again assuming the Bloch--Kato conjecture, 
let $X/\mathbb{C}$ be a variety and $\xi\in H^i(X,\mathbb{Z})$ an element of prime exponent l. Fix some 
points $x_1,\ldots ,x_n\in X$. Then there exists an effective divisor $D\subset X$ such that
$\xi$ restricted to $X-D$ vanishes and $x_j\notin D$ for all $j=1,\ldots ,n$.
\end{coro}

{\em Proof.} This is essentially the same argument as in the proof of Corollary 7.7 of \cite{V2}. \hfill $\Box$

\section{The Milnor--Chow map $\rho$}

\subsection{Higher Chow groups}

S.~Bloch~\cite{B1} defined {\sl higher Chow groups} as a candidate for
motivic cohomology, i.e. an algebraic singular (co)homology. They form a Borel--Moore
homology theory for schemes over a field $k$, which we fix from now on. 
In order to define them we use the algebraic $n$--cube
$$
\square^n=({\mathbb P}^1_k \setminus \{1\})^n. 
$$
The $n$--cube has $2^n$ codimension one faces, defined by $x_i=0$ and 
$x_i=\infty$ for $1 \le i \le n$. An integral subvariety $W \subseteq \square^n$
of codimension $p$ is called admissible if its intersection with all faces 
is again of codimension $p$ or empty.  For each face $F=\{x_i=0\}$ or $F=\{x_i=\infty\}$  
we have a pull--back map $\partial_i^0$ resp. $\partial_i^\infty$ which sends a subvariety 
$W \subseteq \square^n$ to the intersection product of cycles $W \cdot F$ with appropriate 
multiplicities in the sense of Serre's Tor--formula. A total differential is given by
$$
\partial=\sum_{i=1}^n (-1)^{i-1} (\partial_i^0 -\partial_i^\infty).
$$ 
Let $X$ be a quasi--projective variety over $k$ (standard techniques allow 
to extend this definition to equidimensional schemes over $k$ and even, but much harder,  
to schemes over Dedekind rings). The notion of faces, restriction maps, and differentials extends to 
$\square^n_X=X \times_k \square^n$. $Z_c^p(X,n)$ is defined to be the quotient of the group
of admissible cycles  of codimension $p$ in $X \times \square^n$ by the group of degenerate cycles
as defined in \cite{T}, p.180 (where they are denoted by $d^p(X,n)$). 
Let $CH^p(X,n)$ be the $n$--th homology of the complex $Z_c^p(X,\cdot)$ with 
differential $\partial$. 

\subsection{Milnor $K$--theory}

Milnor $K$--theory of a ring $A$  is defined as the quotient 
$$
T(A)/S(A)
$$
of the free graded tensor algebra $T(A)={\mathbb Z} \oplus A^\times \oplus A^\times \otimes A^\times \oplus \cdots$
over the units $A^\times$ of $A$ by the ideal $S(A)$ generated by the degree two relations of the form
$(f,1-f)$ for all $f$ with $f, 1-f \in A^\times$ and $(f,-f)$ for all $f \in A^\times$. 
Note that in the case of fields or (semi--)local
rings with large residue fields the relation $(f,-f)$ follows from the usual Steinberg relation $(f,1-f)$.

\subsection{The map $K^M_n(A)\to CH^n(A,n)$}

Now we consider the special case where $A$ is a localization of an affine $k$--algebra with $k$ an arbitrary ground field.
Denote by $CH^p(A,n)$ the higher Chow groups of ${\rm Spec}(A)$. In particular we have 
the series of abelian groups $CH^n(A,n)$. To any tuple $f=(f_1,\ldots,f_n)$ of elements $f_i \in A^\times$ 
we can associate a map 
$$
f=(f_1,\ldots,f_n): {\rm Spec}(A) \to ({\mathbb P}^1)^n
$$
and hence by restricting to the cube a graph cycle
$$
\Gamma_f={\rm graph}(f_1,\ldots,f_n) \cap \square_A^n.
$$
Since such graph cycles have no boundary, we immediately get a map 
$$
\rho: (A^\times)^n \to CH^n(A,n).
$$
One can show that $\rho$ preserves bilinearity, is skew--commutative, and obeys the Steinberg relations 
$\rho(f,1-f,f_3,\ldots,f_n)=0$ and $\rho(f,-f,f_3,\ldots,f_n)=0$ \cite{PEVSMS}. 
Therefore it descends to a well--defined homomorphism 
$$
\rho: K^M_n(A)\to CH^n(A,n)
$$
for all $n \ge 0$. If $A$ is essentially smooth $CH^*(A,*)$ has a ring structure and $\rho$
becomes a ring homomorphism. In the special case where $A$ is a field $F$ the following result is classical. 

\begin{theo}[Nesterenko/Suslin, Totaro]
$\rho$ is an isomorphism for every field $F$.
\end{theo}

{\em Proof.} Totaro's proof \cite{T} uses cubical higher Chow groups as defined above. 
He shows that any cycle $Z \in CH^n(F,n)$ is equivalent (cobordant) to a norm--cycle
which has all coordinate entries in $F$. This already gives the surjectivity of $\rho$. 
The inverse map $\rho^{-1}$ is defined using the norm as follows: By linearity it is sufficient to define 
$\rho^{-1}$ for $Z$ irreducible. In this case we choose a minimal finite field extension
$L/F$ such that $Z$ corresponds to an $L$--valued point $(z_1,\ldots,z_n)$.
Then $\rho^{-1}(Z)=N_{L/F}(\{z_1,\ldots,z_n\})$ as an element of $K^M_n(F)$, where 
$ N_{L/F}$ is the norm map of Bass and Tate \cite{BT}.
\hfill $\Box$\\

\section{Symbols in general position}

The main result of this section is Proposition 2.8 which in some sense represents the idea that for good extensions
of (semi--)local rings there should be norms of Milnor $K$--groups as in the field case. In fact such norms can be 
constructed by an extension of the methods described below \cite{Ke}.

\subsection{The group $K^t_n(A)$}

Let $A$ be a (semi--)local UFD and $F=Q(A)$ its quotient field. The group $K^t_n(A)$, we are going to define, should be thought of
as the proper Milnor $K$--group of the ring $A[t]_S$, where $S$ denotes the multiplicative system of all
monic polynomials.

\begin{defi}
An $n$--tuple of rational functions 
\[
\left( \frac{p_1}{q_1},\frac{p_2}{q_2},\ldots ,\frac{p_n}{q_n} \right) \in F(t)^n
\]
 with $p_i,q_i\in A[t]$ for $i=1,\ldots n$ is called feasible if
\begin{enumerate}
\item{The highest nonvanishing coefficients of
$p_i,q_i$ are invertible in $A$ for $i=1,\ldots ,n$.}
\item{For every irreducible factor $u$ of $p_i$ or $q_i$ and $v$ of $p_j$ or $q_j$ ($i,j=1,\ldots ,n$, $i\ne j$)
$u$ is either equivalent or coprime to $v$.}
\end{enumerate}   
\end{defi}

Before stating the definition of $K^t_n(A)$ we have to replace ordinary tensor product.

\begin{defi}
Define
\[
\mathcal{T}^t_n(A)=\mathbb{Z}<\{(p_1,\ldots ,p_n)|(p_1,\ldots ,p_n) \text{ feasible, }p_i\in A[t] \text{ irreducible
or unit}\}>/L
\]
\end{defi}
Here $L$ denotes the subgroup generated by elements
\[
(p_1,\ldots ,a p_i,\ldots ,p_n)-(p_1,\ldots ,a ,\ldots ,p_n)-(p_1,\ldots ,p_i,\ldots ,p_n)
\]
with $a\in A^\times$.

By bilinear factorization the element 
\[
(p_1,\ldots ,p_n)\in \mathcal{T}^t_n(A)
\]
is defined for every feasible $n$--tuple with $p_i\in F(t)$.

Now define the subgroup $St\subset \mathcal{T}^t_n(A)$ to be generated by feasible $n$--tuples
\begin{equation}
(p_1,\ldots ,p,1-p,\ldots , p_n)
\end{equation}
and
\begin{equation}
(p_1,\ldots ,p,-p,\ldots ,p_n)
\end{equation}
with $p_i,p\in F(t)$.

\begin{defi}
Define 
$$
K^t_n(A) = \mathcal{T}^t_n(A)/St
$$
\end{defi}

We denote the image of $(p_1,\ldots ,p_n)$ in $K^t_n(A)$ by $\{p_1,\ldots ,p_n\}$.

\subsection{The tame symbol}

Recall that Milnor constructed so called tame symbols
$$
\partial_\pi:K^M_n(F(t)) \longrightarrow K^M_{n-1}(F[t]/(\pi))
$$
for every irreducible $\pi\in F[t]$ \cite{M} -- in fact this construction works for all discrete valuation rings in
contrast to our generalization below.

\begin{prop}[Tame symbol]
For every irreducible, monic polynomial $\pi\in A[t]$ and $n>0$ one has a unique well defined tame symbol
$$
\partial_\pi : K^t_n(A) \longrightarrow K^M_{n-1}(A[t]/(\pi))
$$
which satisfies 
\begin{equation}
\partial_\pi : \{\pi, x_2,\ldots , x_n\} \mapsto \{\bar{x}_2,\ldots ,\bar{x}_n\} 
\end{equation}
for $x_i\in A[t]$ and $x_i$ coprime to $\pi$.

For $\pi=1/t$ there is a similar tame symbol 
$$
\partial_\pi : K^t_n(A) \longrightarrow K^M_{n-1}(A)
$$
which satisfies (3) for $x_i\in A[1/t]$.
\end{prop}

{\em Proof.} Assume $\pi\in A[t]$. Uniqueness is easy to check. 
In order to show existence, introduce according to an idea of Serre a formal 
skew--commutative element $\xi$ with $\xi^2=\xi \{-1\}$ and $\deg(\xi)=1$.
Define a formal map (which is clearly not well defined)
$$
\theta_\pi: \mathcal{T}^t_*(A) \longrightarrow K^M_*(A[t]/(\pi))[\xi]
$$
by 
$$\theta_\pi(u_1 \pi^{i_1} ,\ldots ,u_n \pi^{i_n} ) = (i_1\xi + \{\bar{u}_1\})\cdots (i_n\xi + \{\bar{u}_n\})\; .$$
We define $\partial_\pi$ by taking the (right--)coefficient of $\xi$. This is a well defined homomorphism. 
So what remains to be shown is that
$\partial_\pi$ factors over the Steinberg relations.

Let $x=(\pi^i u,-\pi^i u)$ be feasible, then
\begin{eqnarray*}
\theta_\pi(x)&=&(i\xi +\{\bar{u}\})(i\xi + \{-\bar{u}\}) \\
 & = & i\xi \{-1\} -i\xi \{\bar{u}\} +i\xi \{-\bar{u}\}+\{\bar{u},-\bar{u}\} = 0\; .
\end{eqnarray*}
For $i>0$ and $x=(\pi^i u,1-\pi^i u)$ feasible one has
$$
\theta_\pi(x)= (i\xi +\{\bar{u}\})\{1\} = 0\; .
$$
For $i<0$ and $x=(\pi^{i} u,1-\pi^{i} u)$ feasible one has
\begin{eqnarray*}
\theta_\pi(x) &=& (i\xi + \{\bar{u}\})(i\xi +\{-\bar{u}\})\\
&=& i\xi\{-1\} +i\xi\{-\bar{u}\}-i\xi\{\bar{u}\}+ \{\bar{u},-\bar{u}\} =0\; . 
\end{eqnarray*}
\hfill $\Box$\\

The tame symbols from Proposition 2.4 are compatible with the corresponding symbols of the quotient field of $A$.
This is the content of the next lemma.

\begin{lem}
Fix either an irreducible, monic $\pi\in A[t]$ as above and let $B=A[t]/(\pi)$ and $L=Q(B)$ or set $\pi=1/t$, $B=A$
and $L=F$. The square 
$$
\begin{CD}
K^t_n(A) @>{\partial_\pi}>> K^M_{n-1}(B)\\
@VVV @VVV\\
K^M_n(F(t)) @>>{\partial_\pi}> K^M_{n-1}(L)
\end{CD}
$$
is commutative.

Moreover for monic $\pi\in F[t]$ but $\pi\notin A[t]$ the composition
$$
K^t_n(A)\longrightarrow K^M_n(F(t)) \stackrel{\partial_\pi}{\longrightarrow} K^M_{n-1}(F[t]/(\pi))
$$
vanishes.
\end{lem}

\begin{prop}
If the residue fields of $A$ are infinite the map
$$
\oplus_\pi \partial_\pi: K^t_n(A) \longrightarrow \oplus_\pi K^M_{n-1}(A[t]/(\pi))
$$
is surjective, where the sum is over all monic, irreducible $\pi\in A[t]$.
\end{prop}

In fact the kernel of $\oplus_\pi \partial_\pi$ is precisely $K^M_n(A)$, but this is more difficult to show \cite{Ke}.

{\em Proof.} Consider the filtration $L_d\subset K^t_n(A)$, where $L_d$ is generated by the feasible $(x_1,\ldots ,x_n)$ with
$x_i\in A[t]$ of degree at most $d$. One has to show 
$$
\oplus_{\deg(\pi)=d} \partial_\pi: L_d \longrightarrow  \oplus_{\deg(\pi)=d} K^M_{n-1}(A[t]/(\pi))
$$
is surjective. Fix $\pi$ of degree $d$. For a symbol $\xi=\{\bar{x}_2,\ldots ,\bar{x}_{n}\}\in K^M_{n-1}(A[t]/(\pi))$
we can according to the following sublemma suppose without restriction that 
$\zeta=\{\pi,x_2,\ldots ,x_n\}\in K^t_n(A)$ is well defined assuming $x_i$ to be choosen of degree $d-1$. 
As we have $\partial_{\pi'}(\zeta)=0$ for
$\pi'\ne \pi$, $\deg(\pi')=d$, and $\partial_\pi(\zeta)=\xi$, this proves the proposition.

\begin{sublem}[Gabber]
Given monic $y_1,\ldots ,y_k\in A[t]$ and an arbitrary $x\in A[t]$ coprime to $\pi$ there exists a factorization
$$
x \equiv x' x''   \text{ mod } (\pi) 
$$
such that $x', x''\in A[t]$ have invertible highest coefficients, $\deg(x')=\deg(x'')=d-1$ 
and $x', x''$ are coprime to $y_j$ for $j=1,\ldots ,k$. 
\end{sublem}
{\em Proof.} Using the Chinese remainder theorem and reduction modulo all maximal ideals 
we can assume that $A$ is an infinite field. The moduli space
of factorizations $x \equiv x' x''   \text{ mod } (\pi)$ is a nonempty Zariski open subset of $\mathbb{A}_A^{d}$.
As finite intersections of such subsets contain a rational point, the sublemma is proven.  
\hfill $\Box$\\

\subsection{Norms}

With the notation as above ($B=A[t]/(\pi)$, $F=Q(A)$, $L=Q(B)$) let $i:A\to F$ and 
$j:B\to L$ be the natural embeddings. For the convenience of the reader we recall the construction of
norms $$N_{L/F}:K^M_n(L)\to K^M_n(F) $$ from \cite{BT}. 

Given $\xi\in K^M_n(L)$ choose $\zeta\in K^M_{n+1}(F(t))$ such that $\partial_{\pi'}(\zeta)=0$ for
$\pi'\ne \pi$ and $\partial_\pi(\zeta)=\xi$.
Set $N_{L/F}(\xi)=-\partial_{1/t}(\zeta)$. Kato showed this norm depends only on the isomorphism class of $(L,\xi)$ over
$F$ and is functorial \cite{K}.

\begin{prop} \label{main-prop}
We have $$N_{L/F} (\mathrm{im}(j_*)) \subset \mathrm{im}(i_*)$$
with $i_*, j_*$ the homomorphisms induced on Milnor $K$--groups.
\end{prop} 

{\em Proof.} Given $\xi\in K^M_n(B)$ choose by Lemma 2.6 $\zeta\in K^t_{n+1}(A)$ such that 
$\partial_{\pi'}(\zeta)=0$ for $\pi\ne\pi'\in A[t]$ and $\partial_\pi(\zeta)=\xi$. Set $\xi'=-\partial_{1/t}(\zeta)\in K^M_n(A)$.
It follows from Lemma 2.5 that 
$$   
N_{L/F}(j_*(\xi)) =i_*(\xi')\; . 
$$
\hfill $\Box$

\section{Proof of Theorem \ref{main-theorem}}

Assume that $[Z] \in CH^n(A,n)$ is a higher Chow cycle. We want to construct an element $\xi \in K_n^M(A)$ such that
$\rho(\xi)=[Z]$. 

\begin{lem}
$Z$ is cobordant to a sum of irreducible cycles $Z'$ such that 
\begin{enumerate}
\item $Z' \subset \square_A^n$ does not intersect any face.  
\item With the coordinate functions $t_1,\ldots ,t_n\in \mathcal{O}_{Z'}$ one has $A[t_1,\ldots ,t_i]$ essentially
smooth over $k$ and finite over $A$ for every $1\le i\le n$.
\item $ \mathcal{O}_{Z'}= A[t_1,\ldots ,t_n]$. 
\end{enumerate} 
\end{lem}

{\rm Proof.} This follows immediately from the ``easy moving lemma" of Bloch and Levine~\cite[chap. II, 3.5]{L} and is also 
applied and explained in \cite{PEVSMS}. \hfill $\Box$. 

Without loss of generality we may therefore assume that $Z$ is irreducible and already in good position as in the lemma.
Look at the following diagram:
$$
\begin{CD}
K^M_n(A) @>{i_*}>> K^M_n(F)\\
\rho @VVV \rho @VVV \\
CH^n(A,n) @>>{i_*}> CH^n(F,n) 
\end{CD}
$$
where -- by abuse of notation -- we use the same symbols $\rho$ and $i_*$ for the corresponding maps of rings or fields
and $F$ is the quotient field of $A$.
Since $\rho$ is an isomorphism on the level of fields, we know that there is an element $\tau \in K_n^M(F)$ such that  
$\rho(\tau)=i_*[Z]$. By the description of $\rho^{-1}$ in Totaro's proof of Theorem 1.1, we know that one has 
$\tau=N_{L/F}(\{t_1,\ldots,t_n\})$ where $L$ is the quotient field of $\mathcal{O}_Z$ and $N_{L/F}$ 
is the norm on Milnor $K$--theory of fields. Now look at the consecutive extensions 
$$
A \subset A[t_1] \subset A[t_1,t_2] \subset \ldots \subset A[t_1,\ldots,t_i] \subset \ldots
$$
These rings are all essentially smooth and hence factorial. Each extension is of the type 
$$
A[t_1,\ldots,t_{i+1}]=A[t_1,\ldots,t_i][t]/(\pi_{i+1}).
$$ 
Therefore we may apply Proposition~\ref{main-prop} and conclude that there is an element $\xi$ with $i_*(\xi)=\tau$.  
But the map $i_*:CH^n(A,n)\to CH^n(F,n)$ is injective by \cite{B1} and therefore we have $\rho(\xi)=[Z]$, since $i_*(\rho(\xi))=i_*[Z]$.
\hfill $\Box$.

\end{document}